# The index of a real vector field at an isolated complete intersection singularity

Achim Hennings[1]

**Abstract**: In an unpublished note [H1] we have described a method to obtain a formula for the index of an analytic vector field with (complex) isolated zero on a real analytic hypersurface with (complex) isolated singularity. This formula, like the one of Eisenbud-Levine and Khimshiashvili [AGV] for smooth points, expresses the index by the signature of bilinear forms, which are defined by a local residue symbol (cf. [Ma]). In the complete intersection case, we use a generalized residue symbol, defined for free resolutions in [LJ], in the special case of generalized Koszul complexes to obtain a suitable calculus for the bilinear forms involved.

## 1 Definition of a generalized residue symbol

Let $R, \mathfrak{m}$ be a reduced $d$-dimensional local complex analytic algebra, which represents a complete intersection. Let $\kappa_R = R \cdot \omega$ be the canonical module. There are canonical maps for $l \geq 0$,

(1) $$Ext_R^d(R/\mathfrak{m}^l, \kappa_R) \to H_{\mathfrak{m}}^d(\kappa_R) \xrightarrow{res} \mathbb{C}.$$

Now, we consider 0-dimensional residue class rings $\bar{R}$ of the following kind:

- $a_1, \dots, a_k \in R$ is part of a parameter system of $R$
- $A_1, \dots, A_n \in R^m$, with $n = (d-k) + m - 1$, are the columns of a matrix $A$, having an ideal $I_m(A)$ of $m$ minors which modulo $(a_1, \dots, a_k)$ has the maximal height $n - m + 1 = d - k$.

Then put $\bar{R} := R/(a) + I_m(A)$. This is a zero-dimensional ring, which is resolved over $R/(a)$ by a generalized Koszul complex. A free resolution over $R$ therefore is obtained in the form $L \otimes K \to \bar{R}$ from the generalized Koszul complex of $A$ (cf. [No]),

$$L = L(A)\colon S_{d-k-1}R^m \otimes \Lambda^0 R^n \to \cdots \to S_0 R^m \otimes \Lambda^{d-k-1} R^n \to S_0 R^m \otimes \Lambda^{d-k-1+m} R^n,$$

and the ordinary Koszul complex of $a = (a_1, \dots, a_k)$,

$$K = K(a)\colon \Lambda^0 R^k \to \cdots \to \Lambda^k R^k.$$

The maps within these complexes are defined as follows:

$L(A)\colon \partial = \sum_{i=1}^m x_i^{-1} \otimes y_i \wedge$, resp. $\partial = y_1 \wedge \dots \wedge y_m \wedge$, $y_i := \sum_{j=1}^n A_{ij} e_j$,

$K(a)\colon \partial = \sum_{i=1}^k a_i e_i \wedge$

For our purposes, it is convenient to index the complexes from zero on increasingly, and to choose $\partial \otimes id + (-1)^p id \otimes \partial$ as the map on $L^p \otimes K^q$. Doing so, in the special case $m = 1, n =$

[1] Universität Siegen, Fakultät IV, Hölderlinstraße 3, D-57068 Siegen

$d-k$, the complex $L \otimes K$ is in a simple way isomorphic to the ordinary Koszul complex $K(A_1, \dots, A_n, a_1, \dots, a_k)$.

The resolution $L \otimes K \to \bar{R}$ gives us

$$\kappa_{\bar{R}} = Ext_R^d(\bar{R}, \kappa_R) = H^0(Hom(L \otimes K, \kappa_R)) \cong (S_{d-k-1}\bar{R}^m/(A)_{d-k-1}) \otimes \kappa_R,$$

where $(A) \subseteq SR^m$ is the ideal generated by the columns. We take a surjection $R/\mathfrak{m}^l \to \bar{R}$ with $l \gg 0$ and obtain by (1) the trace map:

$$\tau: \kappa_{\bar{R}} = Ext_R^d(\bar{R}, \kappa_R) \to Ext_R^d(R/\mathfrak{m}^l, \kappa_R) \to \mathbb{C}.$$

For $p\omega \in S_{d-k-1}R^m \otimes \kappa_R$ we now define the generalized residue symbol

$$\operatorname{Res}\begin{bmatrix} p\omega \\ A_1, \dots, A_n, a_1, \dots, a_k \end{bmatrix}$$

as the image of $[p\omega] \in \kappa_{\bar{R}}$ by $\tau$. This definition is a special case of [LJ] for general free resolutions. In our case the symbol has similar properties as the ordinary residue symbol.

By definition, the symbol depends on $p$ only modulo $A_1, \dots, A_n$ and $a_1, \dots, a_k$. The following rules do hold:

(R1) $$\operatorname{Res}\begin{bmatrix} bp\omega \\ A_1, \dots, A_n, ba_1, a_2, \dots, a_k \end{bmatrix} = \operatorname{Res}\begin{bmatrix} p\omega \\ A_1, \dots, A_n, a_1, \dots, a_k \end{bmatrix}$$

for $b \in R$, if the left side is defined.

(R2) $$\operatorname{Res}\begin{bmatrix} bp\omega \\ A_1, \dots, A_{n-1}, bA_n, a_1, \dots, a_k \end{bmatrix} = \operatorname{Res}\begin{bmatrix} p\omega \\ A_1, \dots, A_n, a_1, \dots, a_k \end{bmatrix}$$

for $b \in R$, if the left side is defined.

(R3) $\operatorname{Res}\begin{bmatrix} p\omega \\ A_1, \dots, A_n, a_1, \dots, a_k \end{bmatrix}$ is anti-symmetric in the entries $A_1, \dots, A_n$ and $a_1, \dots, a_k$.

More generally, the symbol satisfies a transformation rule for linear transformations of the $A_i$ resp. $a_j$.

In fact, any homomorphism $\varphi: R^n \to R^n$ induces a homomorphism of complexes, $\varphi_*: L(A) \to L(A^\varphi)$, where $A^\varphi$ is the matrix $A$ with $\varphi$ applied to the rows. The residue map

$$(S_{d-k-1}\bar{R}^m/(A^\varphi)_{d-k-1}) \otimes \kappa_R \to (S_{d-k-1}\bar{R}^m/(A)_{d-k-1}) \otimes \kappa_R$$

therefore represents the multiplication $\det(\varphi): \kappa_{\bar{R}(A^\varphi)} \to \kappa_{\bar{R}(A)}$ induced by $\det(\varphi): \bar{R}(A) \to \bar{R}(A^\varphi)$. Here we have indicated the respective defining matrices. By this remark follows (R2) and part of (R3). A similar reasoning for the complex $K$ proves (R1) and the second part of (R3).

(R4) $$\operatorname{Res}\begin{bmatrix} A_{n+1}p\omega \\ A_1, \dots, A_n, A_{n+1}a_1, a_2, \dots, a_k \end{bmatrix} = \operatorname{Res}\begin{bmatrix} p\omega \\ A_1, \dots, A_n, a_1, \dots, a_k \end{bmatrix}$$

for $A_{n+1} \in S_1R^m$, if the left side is defined.

Proof: We construct a homomorphism

$$L(A_1,\dots,A_n,a_1A_{n+1}) \to L(A_1,\dots,A_n)\otimes K(a_1).$$

Applying $\otimes K(a_2,\dots,a_n)$ we get a homomorphism of the defining complexes for the residue symbols, from which follows the formula. By the decomposition

$$\Lambda^p(R^n\oplus R\cdot e_{n+1}) = \Lambda^{p-1}(R^n)\cdot e_{n+1}\oplus\Lambda^p(R^n)$$

we obtain a representation of $L(A_1,\dots,A_n,a_1A_{n+1})$ as the total complex of the double complex

$$\begin{array}{ccccccc} S_{n-m}\otimes\Lambda^0 & \to\cdots\overset{\partial}{\to} & S_1\otimes\Lambda^{n-m-1} & \overset{\partial}{\to} & S_0\otimes\Lambda^{n-m} & \overset{\alpha}{\to} & S_0\otimes\Lambda^n \\ \uparrow B & & \uparrow(-1)^{n-m-1}B & & \uparrow(-1)^{n-m}B & & \uparrow(-1)^n id\otimes\beta\wedge \\ S_{n-m+1}\otimes\Lambda^0 & \to\cdots\overset{\partial}{\to} & S_2\otimes\Lambda^{n-m-1} & \overset{\partial}{\to} & S_1\otimes\Lambda^{n-m} & \overset{\partial}{\to} & S_0\otimes\Lambda^{n-(m-1)} \end{array}$$

where we have used the notation

$$S := SR^m,\ \Lambda := \Lambda R^n,\ B := a_1A_{n+1}$$

$$(A_{ij}) := A_j,\ y_i := \textstyle\sum_{j=1}^n A_{ij}e_j,\ B = (B_i)$$

$$\alpha := y_1\wedge\dots\wedge y_m,\ \beta := \textstyle\sum_{i=1}^m(-1)^{i-1}y_1\wedge\dots\hat{y}_i\cdots\wedge y_mB_i,\ \partial := \sum_{i=1}^m x_i^{-1}\otimes y_i\wedge$$

and $B$ for the vertical map $\sum_{i=1}^m x_i^{-1}B_i\otimes id$. Analogously, $L(A_1,\dots,A_n)\otimes K(a_1)$ is represented by the double complex:

$$\begin{array}{ccccccc} S_{n-m}\otimes\Lambda^0 & \to\cdots\overset{\partial}{\to} & S_1\otimes\Lambda^{n-m-1} & \overset{\partial}{\to} & S_0\otimes\Lambda^{n-m} & \overset{\alpha}{\to} & S_0\otimes\Lambda^n \\ \uparrow a_1 & & \uparrow(-1)^{n-m-1}a_1 & & \uparrow(-1)^{n-m}a_1 & & \uparrow(-1)^{n-m+1}a_1 \\ S_{n-m}\otimes\Lambda^0 & \to\cdots\overset{\partial}{\to} & S_1\otimes\Lambda^{n-m-1} & \overset{\partial}{\to} & S_0\otimes\Lambda^{n-m} & \overset{\alpha}{\to} & S_0\otimes\Lambda^n \end{array}$$

The homomorphism looked for between the two diagrams is then obtained like this: the first line is mapped identically. The second line is mapped by $A_{n+1}$ (replacing $\pm B$) and by

$$(-1)^{m-1}\gamma,\ \gamma := \textstyle\sum_{i=1}^m(-1)^{i-1}y_1\wedge\dots\hat{y}_i\cdots\wedge y_mA_{i,n+1}$$

(replacing $(-1)^n\beta$). As the dual homomorphism is the identity in total degree $n-m+2 = d-(k-1)$ and multiplication by $A_{n+1}$ in total degree 0, the formula follows.

As a supplement to (R4) we have:

(R5) Let $k = d-1,\ n = (d-k)+m-1 = m,\ \Delta = \det(A_1,\dots,A_n),\ p\in R = S_0R^m$. Then

$$\operatorname{Res}\begin{bmatrix} p\omega \\ \Delta, a_1,\dots,a_k\end{bmatrix} = \operatorname{Res}\begin{bmatrix} p\omega \\ A_1,\dots,A_n,a_1,\dots,a_k\end{bmatrix}.$$

This holds because the defining complex reads in this case

$$L\colon S_0R^m\otimes\Lambda^0R^n \xrightarrow{\Delta} S_0R^m\otimes\Lambda^mR^n.$$

A further conversion of the residue symbol is obtained by transforming the columns of the matrix $A$. For this purpose, it is worthwhile to view the complex $L$ as the dual of the complex with maps

$$\sum_j(\sum_i A_{ij}x_i)\otimes e_j^*\colon S_pR^m\otimes\Lambda^qR^n \to S_{p+1}R^m\otimes\Lambda^{q-1}R^n.$$

For a given isomorphism $\varphi\colon R^m \to R^m$, let $A^\varphi$ be the matrix with columns $\varphi(A_j)$. If we consider $A_j$ as $\sum_i A_{ij}x_i \in S_1R^m$ and $A_j^\varphi$ as $\sum_i(A^\varphi)_{ij}x_i \in S_1R^m$, then $A_j^\varphi = S(\varphi)A_j$. We have a commutative diagram

$$\begin{array}{ccc} S_pR^m\otimes R^n & \xrightarrow{(A_1^\varphi,\dots,A_n^\varphi)} & S_{p+1}R^m \\ \uparrow S_p(\varphi) & & \uparrow S_{p+1}(\varphi) \\ S_pR^m\otimes R^n & \xrightarrow{(A_1,\dots,A_n)} & S_{p+1}R^m \end{array}$$

which can be prolonged to the exterior powers of $R^n$.

Then we have the dual diagrams (for $p+q=d-k-1=n-m$)

$$\begin{array}{ccc} (S_{p+1}R^m\otimes\Lambda^{q-1}R^n)^* & \xrightarrow{\partial_{A^\varphi}} & \left(S_pR^m\otimes\Lambda^qR^n\right)^* \\ \downarrow S_{p+1}(\varphi)^* & & \downarrow S_p(\varphi)^* \\ (S_{p+1}R^m\otimes\Lambda^{q-1}R^n)^* & \xrightarrow{\partial_A} & \left(S_pR^m\otimes\Lambda^qR^n\right)^* \end{array}$$

and at the end

$$\begin{array}{ccc} \Lambda^{n-m}R^n & \xrightarrow{v_1\wedge\dots\wedge v_m\wedge} & \Lambda^nR^n \\ \downarrow id & & \downarrow (\det\varphi)^{-1} \\ \Lambda^{n-m}R^n & \xrightarrow{u_1\wedge\dots\wedge u_m\wedge} & \Lambda^nR^n \end{array}$$

where $u_1,\dots,u_m$ are the rows of $A$ and $v_1,\dots,v_m$ are the rows of $A^\varphi$. Since for calculating $Ext$ we need to dualize once again, we see that the map

$$(S_{d-k-1}\bar{R}^m/(A)_{d-k-1})\otimes\kappa_R \xrightarrow{S(\varphi)} (S_{d-k-1}\bar{R}^m/(A^\varphi)_{d-k-1})\otimes\kappa_R$$

represents the map $\det(\varphi)^{-1}\colon\kappa_{\bar{R}(A)} \to \kappa_{\bar{R}(A^\varphi)}$ induced by $\det(\varphi)^{-1}\colon\bar{R}(A^\varphi)\to\bar{R}(A)$, where $\bar{R}=\bar{R}(A)=\bar{R}(A^\varphi)$ (since $I_m(A)=I_m(A^\varphi)$). We arrive at the formula

$$\operatorname{Res}\begin{bmatrix}\det(\varphi)\,S(\varphi)(p)\omega \\ \varphi(A_1),\dots,\varphi(A_n),a_1,\dots,a_k\end{bmatrix} = \operatorname{Res}\begin{bmatrix}p\omega \\ A_1,\dots,A_n,a_1,\dots,a_k\end{bmatrix}.$$

Furthermore, there is a duality:

(R6) Let $p \in S_i R^m$ and assume for all $q \in S_j R^m$, $i + j = d - k - 1$,

$$\operatorname{Res}\begin{bmatrix} pq\omega \\ A_1, \dots, A_n, a_1, \dots, a_k \end{bmatrix} = 0.$$

Then $p \in (A, a)_i$ if $i > 0$, and $p \in (a) + I_m(A)$ if $i = 0$. Hence,

$$\oplus_{i=0}^{d-k-1} S_i(\bar{R}^m/(A))$$

is self-dual.

Sketch of proof for (R6): Over $\tilde{R} \coloneqq R/(a_1, \dots, a_k)$ we have the exact sequence

$$0 \to S_i R^m \otimes \Lambda^0 R^n \otimes \tilde{R} \to \cdots \to S_0 R^m \otimes \Lambda^i R^n \otimes \tilde{R} \to M \to 0$$

and the sequence

$$0 \to M \to S_0 R^m \otimes \Lambda^{i+m} R^n \otimes \tilde{R} \to \cdots \to S_{n-m-i} R^m \otimes \Lambda^n R^n \otimes \tilde{R} \to 0,$$

which is exact outside the maximal ideal. The first sequence computes $Ext_{\tilde{R}}^i(M, \tilde{R})$ (for $0 < i \le n - m + 1$), the second one computes $H_{\mathfrak{m}}^{n-m-i+1}(M)$ (for $0 < n - m - i < n - m$). By the local duality theorem, these two modules are dual over $\mathbb{C}$ for $0 < i < n - m$.

For $i = 0$, the sequence

$$0 \to S_0 R^m \otimes \Lambda^0 R^n \otimes \tilde{R} \to S_0 R^m \otimes \Lambda^m R^n \otimes \tilde{R} \to S_1 R^m \otimes \Lambda^{m+1} R^n \otimes \tilde{R} \dots$$
$$\to S_{n-m} R^m \otimes \Lambda^n R^n \otimes \tilde{R} \to 0$$

computes the module

$$Ext_{\tilde{R}}^{n-m+1}(\bar{R}, \tilde{R}) \subseteq H_{\mathfrak{m}}^{n-m+1}(\tilde{R}),$$

and this again is dual to $\bar{R} = H_{\mathfrak{m}}^0(\bar{R})$. We do not discuss the compatibility of these dualities with multiplication on $S(\bar{R}^m)/((A) + S_{d-k}(\bar{R}^m))$ here.

Conclusion: In the situation of (R2), $S_i A_n \subseteq S_{i+1}$ and $S_j b \subseteq S_j$ $(i + 1 + j = d - k - 1)$ generate the mutual orthogonal spaces with respect to

$$\operatorname{Res}\begin{bmatrix} M \cdot \omega \\ A_1, \dots, A_{n-1}, A_n b, a_1, \dots, a_k \end{bmatrix},$$

where $M: S_{i+1} \times S_j \to S_{d-k-1}$ is multiplication.

Proof: Assume $p \perp S_j b$. From $pb \perp S_j$ we deduce first

$$pb \equiv q A_n b \bmod (A_1, \dots, A_{n-1}, a_1, \dots, a_k),$$

$$b(p - q A_n) \equiv 0 \bmod (A_1, \dots, A_{n-1}, a_1, \dots, a_k),$$

and then

$$p - qA_n \equiv 0 \bmod (A_1, \dots, A_{n-1}, a_1, \dots, a_k)$$

(since $S_i/(A_1, \dots, A_{n-1}, a_1, \dots, a_k)_i$ for $i \leq d - k - 1$ is a CM module of dimension one by assumption in (R2)).

Assume $p \perp S_i A_n$. From $pA_n \perp S_i$ we infer first

$$pA_n \equiv qA_n b \bmod (A_1, \dots, A_{n-1}, a_1, \dots, a_k),$$

$$A_n(p - qb) \equiv 0 \bmod (A_1, \dots, A_{n-1}, a_1, \dots, a_k),$$

and then

$$p - qb \equiv 0 \bmod (A_1, \dots, A_{n-1}, a_1, \dots, a_k)$$

(since the complex defined by $A_1, \dots, A_n$,

$$S_{j-1}R^m \otimes \Lambda^2 R^n \otimes \tilde{R} \to S_j R^m \otimes \Lambda^1 R^n \otimes \tilde{R} \to S_{j+1} R^m \otimes \Lambda^0 R^n \otimes \tilde{R}$$

is exact for $j < d - k - 1$). A different approach to the proof would be by (R4).

We need a variant of these assertions for real functions. A holomorphic function $f(z)$ (with a symmetric domain) will be called real, if it is invariant by the involution $\sigma(f) = \overline{f(\bar{z})}$, which is an anti-linear ring homomorphism. This condition means for a power series that all coefficients are real. The involution is extended in a multiplicative way to differential forms by $\sigma(df) = d\sigma(f)$ (esp. $\sigma(dz) = dz$).

Now we suppose that the complete intersection ring $R$ is defined by an ideal generated by real power series, so that the involution is defined on $R$ and $\Omega_R^\bullet$. We denote by $R_+ = R^\sigma$ the subalgebra of real elements of $R$. This is the structural algebra of the corresponding real analytic germ, namely the quotient of the real power series ring modulo the same generators.

If $A_1, \dots, A_n$, $a_1, \dots, a_k$ and $\omega$ are real, then the residue

$$\operatorname{Res} \begin{bmatrix} p\omega \\ A_1, \dots, A_n, a_1, \dots, a_k \end{bmatrix}, \text{ for all } p \in S_{d-k-1} R_+^m,$$

is a real number. To see this, we construct $H_{\mathfrak{m}}^d(\kappa_R)$ by means of real parameter systems. Then we have a canonical involution on this module, and the map (1) becomes equivariant.

Under the same hypothesis, the duality of (R6) and the conclusion also hold for the real subspaces $S_i R_+^m$. In fact, the space

$$\oplus_{i=0}^{d-k-1} S_i(\bar{R}^m/(A))$$

is decomposed by $\sigma$ into the eigenspaces $V_+, V_-$ belonging to the eigenvalues $\pm 1$, and we have $iV_+ = V_-$. Therefore, the bilinear form given by the residue symbol stays to be non-degenerate on the subspace $V_+$, and we derive the conclusion as before.

We consider now a deformation of a $d$-dimensional reduced complete intersection $(X,0) \subseteq (\mathbb{C}^N,0)$. Let $X$ be defined by $f = (f_1,\dots,f_c)$ and the deformation $(Y,0) \subseteq (\mathbb{C}^N \times T,0)$ defined by $(f_{t1},\dots,f_{tc})$, where $(T,0) = (\mathbb{C},0)$. Let $\pi: Y \to T$ be the projection and $Y_t$ the fiber. The elements $A_1,\dots,A_n$ and $a_1,\dots,a_k$ may also depend on the parameter $t \in T$. Then the subspace $Z \subseteq Y$ defined by $(a) + I_m(A)$ is finite over $T$. For each fiber we can form the residue sum

$$\mathrm{Res}_{Y_t}\begin{bmatrix} p\omega \\ A_1,\dots,A_n,a_1,\dots,a_k \end{bmatrix}, p \in \Gamma(S_{d-k-1}\mathcal{O}_Y^m).$$

This is a holomorphic function of $t$. For a proof, we choose functions $z_1,\dots,z_N \in \mathcal{O}_{\mathbb{C}^N\times T}$ defining a complete intersection $Z' \supseteq Z$ over $T$, such that there is a comparison map

$$K(z_1,\dots,z_N) \to L\otimes K$$

between the ordinary Koszul complex of $z_1,\dots,z_N$ and the complex $L\otimes K$ used above (and here in a parameter dependent form) in the definition of the generalized residue symbol. Then the image of $[p\omega] \in \kappa_{Z/T} = Ext_Y^d(Z,\kappa_{Y/T})$ in $Ext_Y^d(Z',\kappa_{Y/T}) \subseteq H_{Z'}^d(\kappa_{Y/T})$ is represented by a cocycle with $z_1,\dots,z_N$ in the denominator, and we can use the holomorphic dependence of the ordinary residue symbol on parameters. The module

$$\oplus_{i=0}^{d-k-1}\pi_*S_i(\mathcal{O}_Y^m/((a)+I_m(A))\mathcal{O}_Y^m + (A))$$

is free over $\mathcal{O}_T$ and the residue symbol defines a holomorphic bilinear form thereon.

Once again, we discuss the case where all data are given by real functions (and $X, Y$ are symmetric). Let $X_{\mathbb{R}} \subseteq \mathbb{R}^N$, $Y_{\mathbb{R}} \subseteq \mathbb{R}^N \times \mathbb{R}$ and $T_R$ be the real subspaces of $X, Y$ and $T$. Let $R_+$ be the ring of $\sigma$-invariant functions of $\Gamma(Y,\mathcal{O}_Y)$ and

$$S_iR_+^m \times S_iR_+^m \to S_{2i}R_+^m$$

the multiplication. In the case $2i = d-k-1$ we consider the bilinear form

$$\phi = Res_{Y_t}\begin{bmatrix} M_i\omega \\ A_1,\dots,A_n,a_1,\dots,a_k \end{bmatrix},$$

which is continuous in $t$ and has real values. In particular, the signature is well-defined. In computing the signature, conjugate non-real pairs of zeroes in $Z$ may be discarded. In fact, for such a pair $y_1, y_2 \in Y_t$ we put

$$V_j := S_i(\mathcal{O}_{Y_t,y_j}^m)/(a,I_m(A),A), j = 1,2,$$

and consider the space $V_+ \subseteq V_1 \oplus V_2$ of $\sigma$-invariant elements. Take a basis $v_1,\dots,v_r$ of $V_1$. Then $\sigma v_1,\dots,\sigma v_r$ is a basis of $V_2$, and $v_1+\sigma v_1,\dots,v_r+\sigma v_r$ together with $i(v_1-\sigma v_1),\dots,i(v_r-\sigma v_r)$ constitute a real basis of $V_+$. The matrix of $\phi$ in this basis has block shape

$$\begin{bmatrix} \alpha & \beta \\ \beta & -\alpha \end{bmatrix}.$$

The characteristic polynomial of this matrix is even, hence the signature zero. This is apparent from the characteristic matrix by multiplying all columns of the left half and all rows of the lower half by $-1$ and rearranging the rows and columns.

We shall write

$$\text{Sig Res}_{Y_{\mathbb{R},t}} \begin{bmatrix} M_i\omega \\ A_1, \dots, A_n, a_1, \dots, a_k \end{bmatrix}$$

for the sum of the signatures at all real (or complex) zeroes in $Z \cap Y_t$, where the residue is evaluated on $Y_t$. The bilinear form is well-defined on the free $\mathcal{O}_{T_\mathbb{R}}$-module

$$\pi_*(S_i(\mathcal{O}_Y^m/(a, I_m(A), A)_i)^\sigma$$

of invariant elements, and as such it is non-degenerate on the fibers. Therefore, the signature sum is constant. However, this property may be no longer true, if the non-degeneracy is lost by an additional factor in the numerator of the symbol.

The conversions of the symbol in the next section pursue the aim to represent a given bilinear form solely by non-degenerate ones.

## 2 Calculation of the index

Let $X_\mathbb{C} = f^{-1}(0) \subseteq (\mathbb{C}^n, 0)$ be a complete intersection of dimension $d \geq 1$ with isolated singularity, defined by $f = (f_1, \dots, f_m)$, $n = d + m$. Let $v_\mathbb{C} = (v_1, \dots, v_n)$ be an analytic vector field on $X_\mathbb{C}$ with isolated zero at $x_0 = 0$. We suppose that all component functions of $f$ and $v$ are real, and we intend to study the real subspace $X = (X_\mathbb{C} \cap \mathbb{R}^n) \subseteq (\mathbb{R}^n, 0)$ and the induced vector field $v$ on $X$. If we exclude the trivial case $X = \{0\}$, $X$ is an isolated singularity of dimension $d$.

We briefly recall the definition of the index of a vector field $v$ with isolated zero at an isolated real-analytic singularity $x_0 = 0$. Let $\varepsilon_0 > 0$ small enough, such that the spheres $S_\varepsilon$ for $0 < \varepsilon \leq \varepsilon_0$ are transversal to $X$ and there is no zero $x \neq x_0$ of $v$ with $\|x\| \leq \varepsilon_0$. For $0 < \varepsilon_1 < \varepsilon_2 < \varepsilon_0$ we look at the manifold

$$X_{\varepsilon_1}^{\varepsilon_2} = \{x \in X | \varepsilon_1 \leq \|x\| \leq \varepsilon_2\}$$

with boundary $K_i = X \cap S_{\varepsilon_i}$, $i = 1,2$. Then we take a continuous vector field $u$ on $X_{\varepsilon_1}^{\varepsilon_2}$ with isolated zeroes, which coincides with $v$ on $K_2$ and transversally points to the interior of $X_{\varepsilon_1}^{\varepsilon_2}$ at $K_1$. The definition reads

$$\text{ind}_{x_0} v = 1 + \sum\nolimits_{x \in X_{\varepsilon_1}^{\varepsilon_2}} \text{ind}_x u.$$

This definition satisfies the expected Hopf index formula and deformation properties.

In the following theorem, the prerequisites are like this:

With $z_1, \dots, z_n$ coordinates on $\mathbb{C}^n$, $A = (A_1, \dots, A_n)$ is the Jacobian matrix of $f$ and $\Delta$ is the sub-determinant of $A$ corresponding to $z_1, \dots, z_m$. The form $\omega$ is the unique $d$-form on $X_{\mathbb{C}} \backslash \{0\}$ such that $df_1 \wedge \dots \wedge df_m \wedge \omega = dz_1 \wedge \dots \wedge dz_n$. By $M$ is meant the multiplication map on the algebra of real functions $\mathcal{O}_X$ and $M_k$ is the product $S_k \times S_k \to S_{2k}$ on the symmetric algebra of $\mathcal{O}_X^m$.

$(Y_{\mathbb{C}}, 0) \subseteq (\mathbb{C}^n \times \mathbb{C}, 0)$ is a smoothing of $X_{\mathbb{C}}$ defined by real functions, and $\tilde{\chi}(Y_t)$ is the reduced (i.e. diminished by 1) Euler characteristic of the real part of the fiber.

The coordinates $z_1, \dots, z_n$ must be chosen such that the submodule $(v_1 A_1, \dots, v_n A_n)$ has finite colength (and are real in the above sense).

For $d - 1$ odd, we need that $Y_{\mathbb{C},t} \cap \{z | z_n = 0\}$ is a smoothing of $X_{\mathbb{C}} \cap \{z | z_n = 0\}$.

**Theorem:** There are two cases: If $d - 1 = 2k$ is even, then

$$\mathrm{ind}_0 v = \mathrm{Sig\ Res}_{X_{\mathbb{C}},0} \begin{bmatrix} M\Delta\omega \\ v_{m+1}, \dots, v_n \end{bmatrix} + \mathrm{Sig\ Res}_{X_{\mathbb{C}},0} \begin{bmatrix} M_k v_n \omega \\ A_1, \dots, A_{n-1} \end{bmatrix} - \tilde{\chi}(Y_t).$$

If $d - 1 = 2k + 1$ is odd, then

$$\mathrm{ind}_0 v = \mathrm{Sig\ Res}_{X_{\mathbb{C}},0} \begin{bmatrix} M\Delta\omega \\ v_{m+1}, \dots, v_n \end{bmatrix} - \mathrm{Sig\ Res}_{X_{\mathbb{C}},0} \begin{bmatrix} M_k A_n \omega \\ A_1, \dots, A_{n-1} \end{bmatrix} - \tilde{\chi}(Y_t \cap \{x | x_n = 0\}).$$

## 2.1 Algebraic conversions

For the proof, we first restrict ourselves to such vector fields, which are tangent to the whole map $f: (\mathbb{R}^n, 0) \to (\mathbb{R}^m, 0)$ (to the fibers of the smoothing is enough). Afterwards we deduce the general case from this.

<u>Special case:</u> $v_{\mathbb{C}}$ is the restriction of a vector field in $\Theta_f$, i.e. a vector field tangent to the fibers of $f$, denoted by the same symbol. We say that the original vector field deforms to a vertical vector field.

We consider a one-parameter deformation $Y_{\mathbb{C}} \to T_{\mathbb{C}}$, where $(Y_{\mathbb{C}}, 0) \subseteq (\mathbb{C}^n \times \mathbb{C}, 0)$, of $X_{\mathbb{C}} = Y_{\mathbb{C},0}$ with smooth general fiber $Y_{\mathbb{C},t}$, $t \neq 0$, which is defined by real functions. In particular, we can take for $Y_{\mathbb{C}}$ the inverse image $f^{-1}(T_{\mathbb{C}})$ of the complexification $T_{\mathbb{C}}$ of a suitable real line $T := L \subseteq \mathbb{R}^m$ through the origin, embedded as a graph into $\mathbb{C}^n \times T_{\mathbb{C}}$. The real part $Y \to T$ is a smoothing of $X$, and $v$ is a vertical vector field.

The Euler characteristics in the theorem are independent of the smoothing chosen, as the dimension is odd.

The following residue calculations refer to the complex extension, however we have to sum only over real points as explained above. We put $T' = T \backslash \{0\}$.

By the definition of the index, we can express it in terms of the index sum on the fiber over $t \in T'$,

$$\mathrm{ind}_0 v|Y_0 = \mathrm{ind}\, v|Y_t - \tilde{\chi}(Y_t).$$

Applying the theorem of Eisenbud and Levine and Khimshiashvili ([AGV, vol. I, ch. 5]), we have

$$\mathrm{ind}\, v|Y_t = \mathrm{Sig\ Res}_{Y_t}\begin{bmatrix} M\Delta\omega \\ v_{m+1}, \dots, v_n \end{bmatrix},$$

where $\Delta = \det\left(\frac{\partial f_i}{\partial z_j}; 1 \le i,j \le m\right)$, $\omega = \frac{dz_1\wedge\ldots\wedge dz_n}{df_1\wedge\ldots\wedge df_m}$, i.e. $\omega$ is the relative $d$-form on $Y_{\mathbb{C}}\backslash\{0\}$ with $df_1 \wedge \ldots \wedge df_m \wedge \omega = dz_1 \wedge \ldots \wedge dz_n$.

This formula easily follows at a point, where $z_{m+1}, \dots, z_n$ are coordinates. For an arbitrary choice $z_j, j \in \bar{J}$, with $J = \{j_1 < \cdots < j_m\}, \bar{J} = \{1, \dots, n\}\backslash J$, one can set up the analogous formula with $\Delta_J = \det\left(\frac{\partial f_i}{\partial z_j}; i = 1, \dots, m, j \in J\right)$ replacing $\Delta$, and with the additional sign $sgn\begin{pmatrix} 1, \dots, n \\ J, \bar{J} \end{pmatrix}$. Again, this gives the index for a coordinate system. In general, the value of this formula does not depend on $J$. As an example, we replace $z_{m+1}$ by $z_1$. Writing $A_j = \partial f/\partial z_j$, we have on $Y_t$:

$$v_1A_1 + \cdots + v_{m+1}A_{m+1} + v_{m+2}A_{m+2} + \cdots + v_nA_n = 0, \text{ i.e.}$$

$$v_1A_1 + \cdots + v_{m+1}A_{m+1} \equiv 0\ mod\ (v_{m+2}, \dots, v_n).$$

By Cramer's rule, with $J = \{2, \dots, m+1\}$, we infer the relation $\Delta v_1 \equiv v_{m+1}\Delta_J(-1)^m$, and the assertion follows. We thereby have proved the index formula in general.

Due to the factor $\Delta$ in the numerator, the signature in the above formula may have a jump for $t \to 0$. For that reason, we transform the residue symbol with the rules we have established earlier.

In this proof we write $S_k$ for $\Gamma(S_k(\mathcal{O}_Y^m))$ and $M_k$ for the multiplication $S_k \times S_k \to S_{2k}$, esp. $M = M_0$ for multiplication on $\Gamma(\mathcal{O}_Y)$.

For the columns $A_j = \partial f/\partial z_j$ of the Jacobian matrix, the relation

$$(2) \qquad v_1A_1 + \cdots + v_nA_n \equiv 0\ mod\ (f_{t1}, \dots, f_{tm})$$

is valid (where $f_t = f - t,\ t \in L$). Since $(\Delta)$ and $(v_{m+1})$ are mutually orthogonal spaces $mod\ (\Delta v_{m+1}, v_{m+2}, \dots, v_n)$, we have the following additivity of signatures (cf. [MS]):

$$\mathrm{Sig\ Res}_{Y_t}\begin{bmatrix} M\Delta\omega \\ v_{m+1}, \dots, v_n \end{bmatrix} = s_1 - \mathrm{Sig\ Res}_{Y_t}\begin{bmatrix} Mv_{m+1}\omega \\ \Delta, v_{m+2}, \dots, v_n \end{bmatrix},$$

where

$$s_1 \coloneqq \mathrm{Sig\ Res}_{Y_t}\begin{bmatrix} M\omega \\ \Delta v_{m+1}, v_{m+2}, \dots, v_n \end{bmatrix}.$$

The first term on the right, $s_1$, is continuous at $t = 0$. The second one, by (R5) and (R4), equals

$$-\text{Sig Res}_{Y_t}\begin{bmatrix} Mv_{m+1}\omega \\ A_1,\dots,A_m,v_{m+2},\dots,v_n \end{bmatrix} = -\text{Sig Res}_{Y_t}\begin{bmatrix} MA_{m+1}v_{m+1}\omega \\ A_1,\dots,A_m,A_{m+1}v_{m+2},v_{m+3},\dots,v_n \end{bmatrix},$$

which is by the relation (2) and (R2) also equal to

$$\text{(3)}\quad \text{Sig Res}_{Y_t}\begin{bmatrix} MA_{m+2}v_{m+2}\omega \\ A_1,\dots,A_m,A_{m+1}v_{m+2},v_{m+3},\dots,v_n \end{bmatrix} = \text{Sig Res}_{Y_t}\begin{bmatrix} MA_{m+2}\omega \\ A_1,\dots,A_{m+1},v_{m+3},\dots,v_n \end{bmatrix}$$

By (R6) and the conclusion thereof, $(A_{m+2}) \subseteq S_1$ and $S_1 v_{m+3} \subseteq S_1$ generate the mutual orthogonal spaces with respect to

$$\text{Res}_{Y_t}\begin{bmatrix} M_1\omega \\ A_1,\dots,A_{m+1},A_{m+2}v_{m+3},v_{m+4},\dots,v_n \end{bmatrix}.$$

By (R4) and (R1) this implies

$$\begin{aligned}
\text{Sig Res}_{Y_t}\begin{bmatrix} MA_{m+2}\omega \\ A_1,\dots,A_m,A_{m+1},v_{m+3},\dots,v_n \end{bmatrix} &= \text{Sig Res}_{Y_t}\begin{bmatrix} MA_{m+2}A_{m+2}\omega \\ A_1,\dots,A_{m+1},A_{m+2}v_{m+3},v_{m+4},\dots,v_n \end{bmatrix} \\
&= s_2 - \text{Sig Res}_{Y_t}\begin{bmatrix} M_1v_{m+3}^2\omega \\ A_1,\dots,A_{m+1},A_{m+2}v_{m+3},v_{m+4},\dots,v_n \end{bmatrix} \\
&= s_2 - \text{Sig Res}_{Y_t}\begin{bmatrix} M_1v_{m+3}\omega \\ A_1,\dots,A_{m+1},A_{m+2},v_{m+4},\dots,v_n \end{bmatrix}
\end{aligned}$$

where

$$s_2 = \text{Sig Res}_{Y_t}\begin{bmatrix} M_1\omega \\ A_1,\dots,A_{m+1},A_{m+2}v_{m+3},v_{m+4},\dots,v_n \end{bmatrix}$$

is again continuous. Once more by the relation (2), we can rewrite the second term as

$$\begin{aligned}
\text{(4)}\qquad & -\text{Sig Res}_{Y_t}\begin{bmatrix} M_1v_{m+3}A_{m+3}\omega \\ A_1,\dots,A_{m+1},A_{m+2},A_{m+3}v_{m+4},v_{m+5},\dots,v_n \end{bmatrix} \\
&= \text{Sig Res}_{Y_t}\begin{bmatrix} M_1v_{m+4}A_{m+4}\omega \\ A_1,\dots,A_{m+2},A_{m+3}v_{m+4},v_{m+5},\dots,v_n \end{bmatrix} \\
&= \text{Sig Res}_{Y_t}\begin{bmatrix} M_1A_{m+4}\omega \\ A_1,\dots,A_{m+2},A_{m+3},v_{m+5},\dots,v_n \end{bmatrix}
\end{aligned}$$

Comparing (3) and (4), we have a shift by two in the denominator. For the outcome, we have to distinguish the two cases of the theorem.

First let $n-m-1=2k$ be even. Then by the method we finally get

$$-\text{Sig Res}_{Y_t}\begin{bmatrix} Mv_{m+1}\omega \\ \Delta,v_{m+2},\dots,v_n \end{bmatrix} = s_2+\dots+s_{k+1} - \text{Sig Res}_{Y_t}\begin{bmatrix} M_kv_n\omega \\ A_1,\dots,A_{n-1} \end{bmatrix},$$

where

$$s_i = \text{Sig Res}_{Y_t}\begin{bmatrix} M_i\omega \\ A_1,\dots,A_{m+2i-3},A_{m+2i-2}v_{m+2i-1},v_{m+2i},\dots,v_n \end{bmatrix}$$

is continuous at $t=0$. Combined with the initial step, we obtain

(*) $$\operatorname{Sig\,Res}_{Y_t}\begin{bmatrix} M\Delta\omega \\ v_{m+1},\dots,v_n \end{bmatrix} = s_1+\dots+s_{k+1}-\operatorname{Sig\,Res}_{Y_t}\begin{bmatrix} M_k v_n\omega \\ A_1,\dots,A_{n-1} \end{bmatrix}.$$

The last term is zero for $t\neq 0$. Actually, either some $m$-minor of $A_1,\dots,A_{n-1}$ or some $m$-minor of $A_1,\dots,A_{n-1},A_n$, which contains $A_n$ is non-zero. Then, from the relation (2) by Cramer's rule $v_n\in I_m(A_1,\dots,A_{n-1})+(f_t)$, and the claim follows. (Of course, the term then vanishes also for $t=0$, but this conclusion needs that $v$ is vertical.)

For $n-m-1=2k+1$ odd, we obtain in a similar way

$$-\operatorname{Sig\,Res}_{Y_t}\begin{bmatrix} Mv_{m+1}\omega \\ \Delta,v_{m+2},\dots,v_n \end{bmatrix} = s_2+\dots+s_{k+1}+\operatorname{Sig\,Res}_{Y_t}\begin{bmatrix} M_kA_n\omega \\ A_1,\dots,A_{n-1} \end{bmatrix}$$

and

(**) $$\operatorname{Sig\,Res}_{Y_t}\begin{bmatrix} M\Delta\omega \\ v_{m+1},\dots,v_n \end{bmatrix} = s_1+\dots+s_{k+1}+\operatorname{Sig\,Res}_{Y_t}\begin{bmatrix} M_kA_n\omega \\ A_1,\dots,A_{n-1} \end{bmatrix}.$$

In order to evaluate the last term, we need an intermediate reflection.

## 2.2 Index of a linear form

We consider a smooth complete intersection $X_{\mathbb{C}}$ of dimension $d=n-m$ in $(\mathbb{C}^n,0)$ defined by $f=(f_1,\dots,f_m)$, and we assume that the linear form $g\coloneqq z_n$ has an isolated singularity on $X_{\mathbb{C}}$, where $z_1,\dots,z_n$ are coordinates. If we write $A=\partial f/\partial z_j$ and $A=(A_1,\dots,A_n)$, then $A$ has an invertible $m$-minor, but $(A_1,\dots,A_{n-1})$ has not. Let us suppose that $\Delta\coloneqq\det(A_{n-m+1},\dots,A_n)$ is invertible, so that $z_1,\dots,z_{n-m}$ is a coordinate system on $X_{\mathbb{C}}$.

Furthermore we assume that $f=(f_1,\dots,f_m)$ is real, and we are interested in the linear form $dg=dx_n$ on the real part $X\subseteq X_{\mathbb{C}}$. Our aim is to express the index by a generalized residue symbol.

We denote the multiplication map on the real algebra $R\coloneqq\mathcal{O}_{X,0}$ by $M$. For the standard basis of $R^m$ we use capital letters $X_1,\dots,X_m$ and for the real coordinates $x_1,\dots,x_n$.

Applying Cramer's rule to the system of equations

$$-\left(f_{i1}dz_1+\dots+f_{i,n-m}dz_{n-m}\right)=f_{i,n-m+1}dz_{n-m+1}+\dots+f_{i,n}dz_n,\ i=1,\dots,m,$$

we get

$$\Delta dz_n=-\textstyle\sum_{j=1}^{n-m}\Delta_j\,dz_j,\ \text{with}\ \Delta_j\coloneqq\det(A_{n-m+1},\dots,A_{n-1},A_j).$$

We have to compute

$$(-1)^d\operatorname{ind}_{X,0}(dx_n)=\operatorname{Sig\,Res}\begin{bmatrix} Mdz_1\wedge\dots\wedge dz_{n-m} \\ \Delta_1/\Delta,\dots,\Delta_{n-m}/\Delta \end{bmatrix}=\operatorname{Sig\,Res}\begin{bmatrix} M\Delta^d dz_1\wedge\dots\wedge dz_{n-m} \\ \Delta_1,\dots,\Delta_{n-m} \end{bmatrix}.$$

By (R5) and (R4), noticing $\det(A_{n-m+1},\dots,A_{n-1},A_1)=\Delta_1$, we have

$$\operatorname{Res}\begin{bmatrix} M\Delta^d dz_1 \wedge \dots \wedge dz_{n-m} \\ \Delta_1, \dots, \Delta_{n-m} \end{bmatrix} = \operatorname{Res}\begin{bmatrix} M\Delta^d dz_1 \wedge \dots \wedge dz_{n-m} \\ A_{n-m+1}, \dots, A_{n-1}, A_1, \Delta_2, \dots, \Delta_{n-m} \end{bmatrix} = \operatorname{Res}\begin{bmatrix} M\Delta^d A_n dz_1 \wedge \dots \wedge dz_{n-m} \\ A_{n-m+1}, \dots, A_{n-1}, A_1, A_n\Delta_2, \Delta_3, \dots, \Delta_{n-m} \end{bmatrix}.$$

We contemplate the matrix $(A_{n-m+1}, \dots, A_n, A_2)$. If we duplicate an arbitrary row and evaluate the determinant, we see:

$$\Delta A_2 = \sum_{j=n-m+1}^{n} (-1)^{n-j} \delta_j A_j, \text{ where } \delta_j \coloneqq \det(A_{n-m+1}, \dots \hat{A}_j \dots, A_n, A_2),\ \delta_n = \Delta_2.$$

From the transformation rule and $\Delta A_2 \equiv \Delta_2 A_n \ mod(A_{n-m+1}, \dots, A_{n-1})$ we get:

$$\begin{aligned} &\operatorname{Res}\begin{bmatrix} M\Delta^d A_n dz_1 \wedge \dots \wedge dz_{n-m} \\ A_{n-m+1}, \dots, A_{n-1}, A_1, A_n\Delta_2, \Delta_3, \dots, \Delta_{n-m} \end{bmatrix} \\ &= \operatorname{Res}\begin{bmatrix} M\Delta^d A_n dz_1 \wedge \dots \wedge dz_{n-m} \\ A_{n-m+1}, \dots, A_{n-1}, A_1, \Delta A_2, \Delta_3, \dots, \Delta_{n-m} \end{bmatrix} \\ &= \operatorname{Res}\begin{bmatrix} M\Delta^{d-1} A_n dz_1 \wedge \dots \wedge dz_{n-m} \\ A_{n-m+1}, \dots, A_{n-1}, A_1, A_2, \Delta_3, \dots, \Delta_{n-m} \end{bmatrix} \end{aligned}$$

By another application of (R4), we rewrite this as

$$\operatorname{Res}\begin{bmatrix} M\Delta^{d-1} A_n A_n dz_1 \wedge \dots \wedge dz_{n-m} \\ A_{n-m+1}, \dots, A_{n-1}, A_1, A_2, A_n\Delta_3, \Delta_4, \dots, \Delta_{n-m} \end{bmatrix},$$

and by evaluating $(A_{n-m+1}, \dots, A_n, A_3)$ we again have

$$\Delta A_3 \equiv \Delta_3 A_n \ mod(A_{n-m+1}, \dots, A_{n-1}).$$

Proceeding like this, we reach at

$$\operatorname{Res}\begin{bmatrix} M\Delta^d dz_1 \wedge \dots \wedge dz_{n-m} \\ A_{n-m+1}, \dots, A_{n-1}, A_1, \Delta_2, \dots, \Delta_{n-m} \end{bmatrix} = \operatorname{Res}\begin{bmatrix} M\Delta A_n^{d-1} dz_1 \wedge \dots \wedge dz_{n-m} \\ A_{n-m+1}, \dots, A_{n-1}, A_1, A_2, \dots, A_{n-m} \end{bmatrix} = (-1)^{(n-m)m-(n-m)(m-1)} \operatorname{Res}\begin{bmatrix} M\Delta^2 A_n^{d-1}\omega \\ A_1, \dots, A_{n-1} \end{bmatrix},$$

where $\omega = dz_1 \wedge \dots \wedge dz_n / df_1 \wedge \dots \wedge df_m = (-1)^{(n-m)m} dz_1 \wedge \dots \wedge dz_{n-m}/\Delta$.

We consider the ring $S(R^m)/(A_1, \dots, A_{n-1}, I_m(A_1, \dots, A_{n-1}))$. On applying an automorphism of $R^m$, $(A_{n-m+1}, \dots, A_n)$ becomes the identity matrix. The ring is then isomorphic to

$$R[X_m]/(A_{n,1}, \dots, A_{n,n-m})$$

and multiplication by $X_m \leftrightarrow A_n$, as well as by the unit $\Delta$, defines an isomorphism of homogeneous components.

By this observation, the last formula simplifies to give the result

$$\operatorname{ind}_{X,0}(dx_n) = \operatorname{Sig}\operatorname{Res}_{X_{\mathbb{C}},0}\begin{bmatrix} M_k\omega \\ A_1, \dots, A_{n-1} \end{bmatrix}$$

if $d-1=2k$, and

$$\mathrm{ind}_{X,0}(dx_n) = \mathrm{Sig}\ \mathrm{Res}_{X_{\mathbb{C}},0}\begin{bmatrix} M_k A_n \omega \\ A_1, \dots, A_{n-1}\end{bmatrix}$$

if $d-1=2k+1$.

## 2.3 Representation of the remainder term in (**) by Milnor numbers

In the following, we assume that $Y \subseteq B_\varepsilon$ is a representative in a small open ball $B_\varepsilon$, and we investigate the fibers $Y_t$ and $Y_{t,s} := Y_t \cap g^{-1}(s)$, where $g = x_n$, for $|s| \le \delta$ and $|t| \le \eta$. Here $\delta = \delta(\varepsilon)$ and $\eta = \eta(\varepsilon)$ are chosen small to avoid singularities at the boundary.

We apply the result of the last section to the smooth manifold $Y_t$, $t \neq 0$, and the function $g_t := x_n\colon Y_t \to \mathbb{R}$. We obtain for the index sum

$$\mathrm{Sig}\ \mathrm{Res}_{Y_t}\begin{bmatrix} M_k A_n \omega \\ A_1, \dots, A_{n-1}\end{bmatrix} = \mathrm{ind}_{\,Y_t}(dg_t).$$

The local index of $dg_t$ at a zero $x \in Y_t$ is also the index of the gradient of $g_t$ (in local coordinates) and as such the negative of the reduced Euler characteristic $\tilde{\chi}_-(g_t, x)$ of the (lower) Milnor fiber (cf. [AGV, vol. II, 14.11]). Summing over all critical points, we have

$$\mathrm{ind}_{\,Y_t}(dg_t) = -\textstyle\sum_x \tilde{\chi}_-(g_t, x).$$

With $\delta > 0$ and $|t| \le \eta = \eta(\delta)$ chosen such that $]-\delta, \delta[$ contains all critical values $s_1 < \dots < s_k$ of $g_t^{-1}$, we put $Z_{[-\delta,\delta]} = g_t^{-1}([-\delta, \delta]) \subseteq Y_t$ and $Z_s = g_t^{-1}(s) \subseteq Y_t$. We choose additional points in $[-\delta, \delta]$, $\sigma_0 < s_1 < \sigma_1 < s_2 \dots < s_k < \sigma_k$. Then (by studying neighborhoods of the critical points) we see

$$\chi(Z_{\sigma_{i-1}}) = \chi(Z_{s_i}) + \tilde{\chi}_-(s_i), \chi(Z_{\sigma_i}) = \chi(Z_{s_i}) + \tilde{\chi}_+(s_i)$$

(where $\tilde{\chi}_\pm(s_i)$ is the sum of $\tilde{\chi}_\pm(g_t, x)$ over all critical points in $Z_{s_i}$) and $\chi(Z_{[\sigma_i,\sigma_{i+1}]}) = \chi(Z_{s_{i+1}})$. We deduce from this:

$$\begin{aligned}\chi(Z_{[\sigma_0,\sigma_{i+1}]}) &= \chi(Z_{[\sigma_0,\sigma_i]}) + \chi(Z_{[\sigma_i,\sigma_{i+1}]}) - \chi(Z_{\sigma_i}) \\ &= \chi(Z_{[\sigma_0,\sigma_i]}) + \chi(Z_{s_{i+1}}) - \chi(Z_{s_{i+1}}) - \tilde{\chi}_-(s_{i+1}) \\ &= \chi(Z_{[\sigma_0,\sigma_i]}) - \tilde{\chi}_-(s_{i+1})\end{aligned}$$

$$\begin{aligned}\chi(Z_{[\sigma_0,\sigma_k]}) &= \chi(Z_{[\sigma_0,\sigma_1]}) - \tilde{\chi}_-(s_2) - \dots - \tilde{\chi}_-(s_k) \\ &= \chi(Z_{\sigma_0}) - \tilde{\chi}_-(s_1) - \tilde{\chi}_-(s_2) - \dots - \tilde{\chi}_-(s_k)\end{aligned}$$

To get a representation by another smooth fiber, we use the rearrangement

$$\chi(Z_{\sigma_{i-1}}) = \chi(Z_{s_i}) + \tilde{\chi}_-(s_i) = \chi(Z_{\sigma_i}) - \tilde{\chi}_+(s_i) + \tilde{\chi}_-(s_i)$$

which leads to

$$\chi(Z_{[\sigma_0,\sigma_k]}) = \chi(Z_{\sigma_i}) - \tilde{\chi}_+(s_1) - \cdots - \tilde{\chi}_+(s_i) - \tilde{\chi}_-(s_{i+1}) - \cdots - \tilde{\chi}_-(s_k).$$

In the present case, $\tilde{\chi}_+ = \tilde{\chi}_-$ since the dimension of $Y_t$ is even, and $Z_0$ is smooth since $Y \cap \{x|x_n = 0\}$ has an isolated singularity by hypothesis. Hence

$$\chi(Z_{[-\delta,\delta]}) = \chi(Z_0) - \sum_x \tilde{\chi}_\pm(g_t, x).$$

The left-hand side is also the Euler characteristic of the Milnor fiber $Y_t$. The following lemma is a special case of a general lemma, which enables to describe the Milnor fibration of a complete intersection inductively (cf. [Lo, prop. 5.4]).

**Lemma:** Let $(Y, 0) \subseteq (\mathbb{R}^n, 0)$ be an isolated real analytic singularity. We consider the function $f: Y \to \mathbb{R}$ and the coordinate function $g = x_n: Y \to \mathbb{R}$. We assume that $f, g$ and $f|g^{-1}(0)$ have an isolated singularity. Here $Y \subseteq B_\varepsilon$ is a representative in the open ball $B_\varepsilon$, $Y_t \coloneqq f^{-1}(t), Y_{t,s} \coloneqq Y_t \cap g^{-1}(s)$, $Y_{t,[a,b]} = Y_t \cap g^{-1}([a,b])$. Then, for sufficiently small positive $\varepsilon$, $\delta(\varepsilon)$ and $t = t(\varepsilon, \delta)$ the subset $Y_{t,[-\delta,\delta]}$ is a deformation retract of the Milnor fiber $Y_t$.

Proof: We deform $Y_t$ to $Y_{t,[-\delta,\delta]}$. By the curve-selection lemma [Lo, 2.1], on $Y_0 \backslash g^{-1}(0)$ the differential of the function $g^2(x) = x_n^2$ is nowhere a non-positive multiple of the differential of $\rho(x) \coloneqq \|x\|^2$, if $\varepsilon$ is sufficiently small. (In fact, otherwise there would be an analytic curve with limit 0, along which $g^2$ does not tend to zero.) Then one can choose $\delta$ and $t = t(\delta)$ such that the stated property remains true on $Y_t \backslash g^{-1}(]-\delta, \delta[)$ and $Y_{t,s}$ is smooth for $|s| \geq \delta$. It is therefore possible to choose a vector field on $Y_t \backslash g^{-1}(]-\delta, \delta[)$, tangent to $Y_t$, which takes on the value $-1$ on $dg^2$ and is not positive on $d\rho$. Then we deform $Y_t$ to $Y_{t,[-\delta,\delta]}$ along the phase curves.

In summary, we have the formula for the remainder term:

$$\text{Sig Res}_{Y_t} \begin{bmatrix} M_k A_n \omega \\ A_1, \ldots, A_{n-1} \end{bmatrix} = \tilde{\chi}(Y_t) - \tilde{\chi}(Y_t \cap \{x|x_n = 0\})$$

*Proof of the special case of the theorem*

We can now confirm the formula in the theorem for a vertical vector field (i.e. tangent to the fibers of the deformation), in particular for nonsingular $X_\mathbb{C}$. To do so, we evaluate the equations (*) and (**) for $t = 0$ and $t \neq 0$, where the (constant) $s_i$ are the same in both cases, and replace the remainder term for $t \neq 0$.

## 2.4 The general case

Here $v_\mathbb{C} \in \Theta_{X_\mathbb{C}}$ need not be deformable to a vertical field.

To extend the index formula to the general case, we need a deformable vector field for comparison. The existence of such a field with isolated zero is affirmed by the following lemma. We restrict to local CM rings $R, \mathfrak{m}$ (of dimension $d = n - m$) as to have equality between height and codimension.

**Lemma 1:** Let $A = (A_1, \dots, A_n)$ be an $m \times n$-matrix without units with the (generic) property

$$ht\, I_m(A_1, \dots, A_k) \geq k - m + 1,\, k = m + 1, \dots, n - 1.$$

Then there is a determinantal relation of the columns $(a_1, \dots, a_n)$ with $ht(a_1, \dots, a_n) \geq n - m$. For any $r > 0$, the relation can be chosen in $\mathfrak{m}^r R^n$.

By a determinantal relation we mean a linear combination of relations arising from the maximal minors of $m + 1$ selected columns (by duplicating an arbitrary row and taking the determinant).

Proof: We show by induction on $m + 1 \leq k \leq n$, that there is a determinantal relation of $A_1, \dots, A_k$ with height $\geq k - m$.

$k = m + 1$: Let $\Delta_i = (-1)^{i+1}\det(A_1, \dots, \hat{A}_i, \dots, A_{m+1})$. Then $\sum_{i=1}^{m+1} \Delta_i A_i = 0$ and $ht(\Delta_1, \dots, \Delta_{m+1}) \geq (m + 1) - m + 1 > k - m$.

$k \to k + 1$: There is a determinantal relation $(a_1, \dots, a_k)$ of $A_1, \dots, A_k$ with height $ht(a_1, \dots, a_k) \geq k - m$. Since $ht\, I_m(A_1, \dots, A_k) \geq k - m + 1$, there is a linear combination $c$ of $m$-minors of $A_1, \dots, A_k$ with $ht(c, a_1, \dots, a_k) \geq k - m + 1$. Hence, there is a determinantal relation $b_1, \dots, b_k, b_{k+1} = c$ of $A_1, \dots, A_k, A_{k+1}$ with height $ht(b_{k+1}, a_1, \dots, a_k) \geq k - m + 1$. Then also the relation

$$(a_1 + tb_1, \dots, a_k + tb_k, tb_{k+1})$$

for general $t$ has at least this height, namely

$$ht(b_{k+1}, a_1 + tb_1, \dots, a_k + tb_k) \geq k - m + 1$$

by semi-continuity of fiber dimension. The claim follows.

To satisfy the additional condition, we can multiply the relation $(\Delta_1, \dots, \Delta_{m+1})$ by a suitable power $\Delta_i^l$ and maintain height one. Similarly, we can multiply the relation $(b_1, \dots, b_k, b_{k+1})$ by a power of $b_{k+1}$ without affecting the height condition of the result.

**Lemma 2:** Let $A = (A_1, \dots, A_n)$ be an $m \times n$-matrix without units with $ht\, I_m(A_1, \dots, A_n) \geq n - m$. Then, by a general column transformation, we can fulfill the hypothesis of the lemma.

Proof: By induction assume that $ht\, I_m(A_1, \dots, A_{k+1}) \geq k - m + 1$ for some $m + 1 \leq k \leq n - 1$. With $c = d - (k - m + 1)$ we find elements $y_1, .., y_c$ such that $\bar{R} = R/(y_1, .., y_c)$ is of dimension $k - m + 1$ and $I_m(A_1, \dots, A_{k+1})\bar{R}$ is zero-dimensional. Over $\bar{R}$ we can arrange that $(A_1, \dots, A_k)$ is a reduction of $(A_1, \dots, A_{k+1})$. Then $I_m(A_1, \dots, A_k)\bar{R}$ is also zero-dimensional, and we have shown the induction step $ht\, I_m(A_1, \dots, A_k) \geq k - m + 1$.

*Proof of the formula in the general case*

Claim 1: For any $k \geq 0$ there is a real deformable vector field $u \in \mathfrak{m}^k \Theta_{X_{\mathbb{C}}}$ with isolated zero on $X_{\mathbb{C}}$.

This is the content of the two lemmas.

Claim 2: For sufficiently large $k > 0$ there is equality $\mathrm{ind}_{X,0}(v+u) = \mathrm{ind}_{X,0}(v)$ for all real $u \in \mathfrak{m}^k \Theta_{X_\mathbb{C}}$ with isolated zero on $X_\mathbb{C}$.

Proof. There is no branching of the zeroes for the homotopy $v + tu$, $t \in [0,1]$, if the norms of the coefficients satisfy $\|u(x)\|^2/\|v(x)\|^2 \to 0$ as $x \to 0$. This is in fact true, if $(u_1, \dots, u_n)\mathcal{O}_{X_\mathbb{C},0} \subseteq \mathfrak{m}(v_1, \dots, v_n)\mathcal{O}_{X_\mathbb{C},0}$.

Claim 3: For sufficiently large $k > 0$ the following holds: For all real $u \in \mathfrak{m}^k \Theta_{X_\mathbb{C}}$ with isolated zero on $X_\mathbb{C}$ and $w(t) = v + tu$, $t \in [0,1]$, the summands $s_i(w(t))$ in (*) and (**) are independent of $t \in [0,1]$.

Proof: This follows from the fact, that for sufficiently large $k \in \mathbb{N}$ the ideals of $\mathcal{O}_{X_\mathbb{C}}$,

$$I_m(A_1, \dots, A_l) + (w_{l+2}, \dots, w_n),\ (\Delta, w_{m+2}, \dots, w_n),\ (w_{m+1}, \dots, w_n),$$

have no branching of the zeroes and therefore the residues are concentrated at the base point 0.

*Conclusion of the proof*

We choose a real deformable vector field $u \in \mathfrak{m}^k \Theta_{X_\mathbb{C}}$ with isolated zero on $X_\mathbb{C}$ for some sufficiently large $k > 0$. We define the vector field $v' = u + \tau v$ for $\tau$ near zero. This may have, besides 0, new zeroes at smooth points of $X$. (The non-real zeroes are not needed.) Then we have continuity (in $\tau$) of:

- The total index $\mathrm{ind}_X(v'(\tau))$.
- The total signatures (expressible by the continuous $s_i(v'(\tau))$)

$$\mathrm{Sig\,Res}_X \begin{bmatrix} M v'_{m+1}\omega \\ \Delta, v'_{m+2}, \dots, v'_n \end{bmatrix} + \mathrm{Sig\,Res}_X \begin{bmatrix} M_k v'_n \omega \\ A_1, \dots, A_{n-1} \end{bmatrix},$$

  resp.

$$\mathrm{Sig\,Res}_X \begin{bmatrix} M v'_{m+1}\omega \\ \Delta, v'_{m+2}, \dots, v'_n \end{bmatrix} - \mathrm{Sig\,Res}_X \begin{bmatrix} M_k A_n \omega \\ A_1, \dots, A_{n-1} \end{bmatrix}.$$

From the smooth case ($\tau \neq 0$ and $x \in X$, $x \neq 0$) and the deformable case ($\tau = 0$) we obtain the validity of the formula for $v' = u + \tau v$, $|\tau| \leq \varepsilon$. Obviously, this implies the validity for $v + u/\varepsilon$. By claim 2 and claim 3 we obtain the validity of the formula for $v$.

## 3 Supplementary remarks

### 3.1 Index of 1-forms

By the reasoning of section 2.2, we can also calculate the index of a general 1-form if the dimension of the space is odd.

We begin as earlier with a nonsingular space, which is given as a complete intersection (in a neighborhood of the base point 0) by real functions,

$$X_\mathbb{C} = f^{-1}(0) \subseteq \mathbb{C}^n, f = (f_1, \dots, f_m), d = n - m,$$

and we are interested in the real subspace $X$. With $A_j := \partial f/\partial z_j$ we form the Jacobi matrix $A = (A_1, \dots, A_n)$. We suppose that $\Delta := \det(A_{n-m+1}, \dots, A_n)$ is invertible, hence $z_1, \dots, z_{n-m}$ a local and system of coordinates.

We consider a real 1-form with isolated zero

$$u = u_1 dz_1 + \cdots + u_n dz_n,$$

and we first express the index of it in local coordinates. Thereafter, we again convert the expression so that the limit to a singular fiber can be recognized.

From the equation $\sum_{j=1}^{n} A_j dz_j = 0$ on $X_{\mathbb{C}}$ we obtain

$$\Delta dz_k = -\sum_{j=1}^{n-m} \Delta_{kj}\, dz_j,\ k = n-m+1, \dots, n,$$

with $\Delta_{kj} = \det(A_{n-m+1}, \dots, A_{k-1}, A_j, A_{k+1}, \dots, A_n)$ (the column with index $k$ is replaced by the one with index $j$). Then

$$u = -\sum_{j=1}^{n-m} \tilde{u}_j\, dz_j \quad \text{with} \quad \tilde{u}_j := u_j - \sum_{k=n-m+1}^{n} \frac{\Delta_{kj}}{\Delta} u_k,$$

and we have

$$ind_{X,0} u = \mathrm{Sig}\ \mathrm{Res}_{X_{\mathbb{C}},0} \begin{bmatrix} M dz_1 \wedge \dots \wedge dz_{n-m} \\ \tilde{u}_1, \dots, \tilde{u}_{n-m} \end{bmatrix}.$$

This formula rests on local coordinates and is restricted to the smooth case. If $X_{\mathbb{C}}$ is a complete intersection with isolated singularity (defined by real functions), we have to take again a deformation $Y \to T$ of $X$ with smooth generic fiber $Y_t$ (and complex extension $Y_{\mathbb{C}}$) and we need a formula which is continuous in $t$.

To find such a formula, we wish to convert the coordinate formula for the smooth $X_{\mathbb{C}}$. We consider the extended columns

$$B_j := \begin{pmatrix} A_j \\ u_j \end{pmatrix}$$

as elements of the free module with standard base $X_1, \dots, X_{m+1}$. Since in the complex case the index is given by the colength of $I_{m+1}(B_1, \dots, B_n)$, it is natural, in the case $n - m - 1 = 2l$, to look at the signature of

$$\mathrm{Res}_{X_{\mathbb{C}},0} \begin{bmatrix} M_l \omega \\ B_1, \dots, B_n \end{bmatrix},$$

$$\omega := dz_1 \wedge \dots \wedge dz_n / df_1 \wedge \dots \wedge df_m.$$

Since $A_{n-m+1}, \dots, A_n$is a basis, we have the representation

$$A_j = \sum_{k=n-m+1}^{n} a_{jk}\, A_k \text{ with } a_{jk} = \Delta_{kj}/\Delta.$$

From this follows

$$\tilde{B}_j := \tilde{u}_j X_{m+1} = B_j - \sum_{k=n-m+1}^{n} \frac{\Delta_{kj}}{\Delta} B_k.$$

By the transformation formula (cf. R3) the above residue symbol becomes

$$\mathrm{Res}_{X_{\mathbb{C}},0}\begin{bmatrix} M_l\omega \\ \tilde{B}_1,\dots,\tilde{B}_{n-m},B_{n-m+1},\dots,B_n\end{bmatrix}.$$

The signature is the same as of

$$\mathrm{Res}_{X_{\mathbb{C}},0}\begin{bmatrix} MX_{m+1}^{2l}\omega \\ \tilde{u}_1X_{m+1},\dots,\tilde{u}_{n-m}X_{m+1},B_{n-m+1},\dots,B_n\end{bmatrix} =$$
$$(-1)^{(n-m)m}\mathrm{Res}_{X_{\mathbb{C}},0}\begin{bmatrix} MX_{m+1}^{2l}\omega \\ B_{n-m+1},\dots,B_n,\tilde{u}_1X_{m+1},\dots,\tilde{u}_{n-m}X_{m+1}\end{bmatrix}.$$

By (R4) the last residue symbol can be simplified as

$$\mathrm{Res}_{X_{\mathbb{C}},0}\begin{bmatrix} M\omega \\ B_{n-m+1},\dots,B_n,\tilde{u}_1X_{m+1},\tilde{u}_2,\dots,\tilde{u}_{n-m}\end{bmatrix}.$$

To the leftmost $m+1$ entries of the denominator, we apply (R5) to get

$$\mathrm{Res}_{X_{\mathbb{C}},0}\begin{bmatrix} M\omega \\ \tilde{u}_1\Delta,\tilde{u}_2,\dots,\tilde{u}_{n-m}\end{bmatrix}.$$

Because of $\omega = (-1)^{(n-m)m}dz_1 \wedge \dots \wedge dz_{n-m}/\Delta$, this is

$$(-1)^{(n-m)m}\mathrm{Res}_{X_{\mathbb{C}},0}\begin{bmatrix} M\Delta^{-2}dz_1 \wedge \dots \wedge dz_{n-m} \\ \tilde{u}_1,\tilde{u}_2,\dots,\tilde{u}_{n-m}\end{bmatrix}.$$

The quadratic factor $\Delta^{-2}$ has no influence on the signature, thus we obtain (for $n-m-1=2l$ even)

$$\mathrm{ind}_{X,0}u = \mathrm{Sig}\,\mathrm{Res}_{X_{\mathbb{C}},0}\begin{bmatrix} M_l\omega \\ B_1,\dots,B_n\end{bmatrix}.$$

This formula, applied to the fibers $Y_t$, is continuous at $t=0$. Therefore, the formula gives the index for isolated singularities $(X_{\mathbb{C}},0)$ if we add the correction term $-\tilde{\chi}(Y_t)$, $t \neq 0$.

Remarks: 1) This formula is an analogue of the formula for the complex index, which is (leaving out the summand $\pm\mu$)

$$(-1)^d \dim \mathcal{O}_{X_{\mathbb{C}},0}/I_{m+1}(B_1,\dots,B_n) = (-1)^d \dim S_{d-1}(\mathcal{O}_{X_{\mathbb{C}},0}^{m+1}/(B_1,\dots,B_n)).$$

2) In the case $n-m-1=2l+1$ odd, for the smooth fibers can be derived a similar formula as above, but with an additional factor $X_{m+1}$ in the numerator. This may cause a jump at $t=0$.

### 3.2 Index of a complex vector field

The method employed here is also suitable to obtain a formula for the index of a holomorphic vector field on a complex complete intersection.

We consider an isolated singularity $(X,0)$ of dimension $d=n-m$ in $\mathbb{C}^n$, defined by $f=(f_1,\dots,f_m)$, and a one-parameter smoothing $Y \to T$ with fibers $Y_t$ of it. To begin with, we consider a vector field $v=\sum_{j=1}^n v_j\,\partial/\partial z_j$ with isolated zero, which is deformable. With $A_j = \frac{\partial f}{\partial z_j}$ we assume the coordinates chosen such that $(v_1A_1,\dots,v_nA_n)$ is of finite codimension. Similar to the real case, we have for $t \neq 0$

$$\mathrm{ind}_{X,0}v = L(\Delta\mathcal{O}_{Y_t}/(v_{m+1},\dots,v_n)) - \tilde{\chi}(Y_t),$$

where $\Delta = \det(A_1,\dots,A_m)$ and $L$ denotes the total length (dimension). In the following, we investigate a fixed fiber $Y_t$ ($t \neq 0$ or $t = 0$) and we shorten $\mathcal{O} := \mathcal{O}_{Y_t}$. For $t = 0$ we do not need that $v$ is deformable.

By nesting of vector spaces we have

$$L(\Delta\mathcal{O}/(v_{m+1},\dots,v_n)) = L(\mathcal{O}/(v_{m+1},\dots,v_n)) - L(\mathcal{O}/(\Delta, v_{m+2},\dots,v_n)) + \\ L(v_{m+1}\,\mathcal{O}/(\Delta, v_{m+2},\dots,v_n)).$$

The first two terms on the right are already part of the expected formula. We write them as multiplicities

$$e(v_{m+1},\dots,v_n) = e(A_1,\dots,A_{m-1},v_{m+1},\dots,v_n),\ e(A_1,\dots,A_m,v_{m+2},\dots,v_n).$$

Such mixed multiplicities are treated in general in [H2], in the present CM-case they are simply the colength of the module generated, modulo the given parameter elements (if there are $< d$ of them). For the third term, we show the injectivity of

$$\mathcal{O}/(\Delta, v_{m+3},\dots,v_n) \xrightarrow{A_{m+1}} \mathcal{O}^m/(A_1,\dots,A_m,v_{m+3},\dots,v_n)).$$

In fact, from $aA_{m+1} \equiv \sum_{j=1}^{m} a_jA_j \ mod(v_{m+3},\dots,v_n)$ we get $a_j\Delta \equiv a\Delta_j$ with $\Delta_j := \det(A_1,\dots,A_{j-1},A_{m+1},A_{j+1},\dots,A_m)$ by Cramer's rule. Since $\mathcal{O}/(\Delta, v_{m+3},\dots,v_n)$ is a one-dimensional CM-ring, and $(\Delta_1,\dots,\Delta_m)$ generates a zero-dimensional ideal, we conclude $a \in (\Delta, v_{m+3},\dots,v_n)$. Now we obtain for the third term:

$$L(v_{m+1}\,\mathcal{O}/(\Delta, v_{m+2},\dots,v_n)) = L(v_{m+1}A_{m+1}\mathcal{O}/(A_1,\dots,A_m,A_{m+1}v_{m+2},v_{m+3},\dots,v_n)) = \\ L(v_{m+2}A_{m+2}\mathcal{O}/(A_1,\dots,A_m,A_{m+1}v_{m+2},v_{m+3},\dots,v_n)) = \\ L(A_{m+2}\mathcal{O}/(A_1,\dots,A_m,A_{m+1},v_{m+3},\dots,v_n)),$$

where we have used in the last step, that $v_{m+2}$ is not a zero-divisor in $\mathcal{O}^m/(A_1,\dots,A_m,v_{m+3},\dots,v_n)$. Again, by nesting, we decompose this as

$$L(\mathcal{O}^m/(A_1,\dots,A_{m+1},v_{m+3},\dots,v_n)) - L(\mathcal{O}^m/(A_1,\dots,A_{m+1},A_{m+2},v_{m+4},\dots,v_n)) + \\ L(v_{m+3}\mathcal{O}^m/(A_1,\dots,A_{m+1},A_{m+2},v_{m+4},\dots,v_n)),$$

the first two terms, $e(A_1,\dots,A_{m+1},v_{m+3},\dots,v_n)$, $e(A_1,\dots,A_{m+2},v_{m+4},\dots,v_n)$, being part of the desired formula. For the third term we now need the injectivity of

$$\mathcal{O}^m/(A_1,\dots,A_{m+2},v_{m+5},\dots,v_n) \xrightarrow{A_{m+3}} S_2\mathcal{O}^m/(A_1,\dots,A_{m+2},v_{m+5},\dots,v_n).$$

It is obtained from the Koszul complex of $A_1,\dots,A_{m+2},A_{m+3}$

$$0 \to S_1 \otimes \Lambda^{3+m} \to S_0 \otimes \Lambda^{2+m} \to S_0 \otimes \Lambda^2 \to S_1 \otimes \Lambda^1 \to S_2 \otimes \Lambda^0 \ \ mod(v_{m+5},\dots,v_n).$$

Applying the injective map, we convert the third term into

$$L(v_{m+3}A_{m+3}S_1\mathcal{O}^m/(A_1,\dots,A_{m+2},A_{m+3}v_{m+4},v_{m+5},\dots,v_n)) = \\ L(v_{m+4}A_{m+4}S_1\mathcal{O}^m/(A_1,\dots,A_{m+2},A_{m+3}v_{m+4},v_{m+5},\dots,v_n)) = \\ L(A_{m+4}S_1\mathcal{O}^m/(A_1,\dots,A_{m+2},A_{m+3},v_{m+5},\dots,v_n)).$$

Repeating the previous steps, we decompose this term by nesting. We get the length of two quotient modules of $S_2\mathcal{O}^m$, which we keep and write as multiplicities since they are Euler characteristics, and a third term to be transformed further. In this way, we obtain an alternating sum of multiplicities and a remainder term. For $d-1=2k$ even, this remainder is $L(v_n S_{d-1}\mathcal{O}^m/(A_1,\dots,A_{n-1}))$, for $d-1=2k+1$ odd it is $L(A_n S_{d-2}\mathcal{O}^m/(A_1,\dots,A_{n-1}))$. The first term vanishes for $t\neq 0$, in view of the relation $v_1A_1+\cdots+v_nA_n=0$, and the second one is for $t\neq 0$ equal to $L(S_{d-1}\mathcal{O}^m/(A_1,\dots,A_{n-1}))$, because $A_n$ induces an isomorphism at the zeroes of $I_m(A_1,\dots,A_{n-1})$ (cf. section 2.2). We note, that this constant multiplicity $e(A_1,\dots,A_{n-1})$ is (by [Lo]) also the sum of Milnor numbers $\mu(f)+\mu(f|z_n=0)$.

Thus, for the fiber $t\neq 0$ and a deformable vector field, we have obtained the formula

$$\mathrm{ind}_{Y_t}v=e(v_{m+1},\dots,v_n)-e(A_1,\dots,A_m,v_{m+2},\dots,v_n)+\cdots+(-1)^d e(A_1,\dots,A_{n-1}).$$

The multiplicities are continuous at $t=0$. For the special fiber the index differs by $\tilde{\chi}(Y_t)$, $t\neq 0$. So the formula reads

$$\mathrm{ind}_{X,0}v=e(v_{m+1},\dots,v_n)-e(A_1,\dots,A_m,v_{m+2},\dots,v_n)+\cdots+(-1)^d e(A_1,\dots,A_{n-1})-\tilde{\chi}(Y_t).$$

We can use the same arguments as in the real case to extend this formula to non-deformable vector fields.